\documentclass[12pt]{amsart}
\usepackage{amssymb}
\begin{document}
%%%%%%%%%%%%%%%%%%%% Text italic %%%%%%%%%%%%%%%%%%%%%%%%%%%%
\theoremstyle{plain}
\newtheorem{thm}{Theorem}[section]
\newtheorem{theorem}[thm]{Theorem}
\newtheorem{lemma}[thm]{Lemma}
\newtheorem{corollary}[thm]{Corollary}
\newtheorem{proposition}[thm]{Proposition}
\newtheorem{addendum}[thm]{Addendum}
\newtheorem{variant}[thm]{Variant}
%%%%%%%%%%%%%%%%%%%% Text roman %%%%%%%%%%%%%%%%%%%%%%%%%%%%%
\theoremstyle{definition}
\newtheorem{notations}[thm]{Notations}
\newtheorem{question}[thm]{Question}
\newtheorem{problem}[thm]{Problem}
\newtheorem{remark}[thm]{Remark}
\newtheorem{remarks}[thm]{Remarks}
\newtheorem{definition}[thm]{Definition}
\newtheorem{claim}[thm]{Claim}
\newtheorem{assumption}[thm]{Assumption}
\newtheorem{assumptions}[thm]{Assumptions}
\newtheorem{properties}[thm]{Properties}
\newtheorem{example}[thm]{Example}
\numberwithin{equation}{thm}
%%%%%%%% Diagram macros, etc. %%%%%%%%%%%%%%%%%%%%%%%%%%%%%%%
\catcode`\@=11
% General macros
\def\opn#1#2{\def#1{\mathop{\kern0pt\fam0#2}\nolimits}}
\def\bold#1{{\bf #1}}%
\def\underrightarrow{\mathpalette\underrightarrow@}
\def\underrightarrow@#1#2{\vtop{\ialign{$##$\cr
 \hfil#1#2\hfil\cr\noalign{\nointerlineskip}%
 #1{-}\mkern-6mu\cleaders\hbox{$#1\mkern-2mu{-}\mkern-2mu$}\hfill
 \mkern-6mu{\to}\cr}}}
\let\underarrow\underrightarrow
\def\underleftarrow{\mathpalette\underleftarrow@}
\def\underleftarrow@#1#2{\vtop{\ialign{$##$\cr
 \hfil#1#2\hfil\cr\noalign{\nointerlineskip}#1{\leftarrow}\mkern-6mu
 \cleaders\hbox{$#1\mkern-2mu{-}\mkern-2mu$}\hfill
 \mkern-6mu{-}\cr}}}
% Rectangular Commutative diagrams
\let\amp@rs@nd@\relax
\newdimen\ex@
\ex@.2326ex
\newdimen\bigaw@
\newdimen\minaw@
\minaw@16.08739\ex@
\newdimen\minCDaw@
\minCDaw@2.5pc
\newif\ifCD@
\def\minCDarrowwidth#1{\minCDaw@#1}
\newenvironment{CD}{\@CD}{\@endCD}
\def\@CD{\def\A##1A##2A{\llap{$\vcenter{\hbox
 {$\scriptstyle##1$}}$}\Big\uparrow\rlap{$\vcenter{\hbox{%
$\scriptstyle##2$}}$}&&}%
\def\V##1V##2V{\llap{$\vcenter{\hbox
 {$\scriptstyle##1$}}$}\Big\downarrow\rlap{$\vcenter{\hbox{%
$\scriptstyle##2$}}$}&&}%
\def\={&\hskip.5em\mathrel
 {\vbox{\hrule width\minCDaw@\vskip3\ex@\hrule width
 \minCDaw@}}\hskip.5em&}%
\def\verteq{\Big\Vert&&}%
\def\noarr{&&}%
\def\vspace##1{\noalign{\vskip##1\relax}}\relax\iffalse{%
\fi\let\amp@rs@nd@&\iffalse}\fi
 \CD@true\vcenter\bgroup\relax\iffalse{%
\fi\let\\=\cr\iffalse}\fi\tabskip\z@skip\baselineskip20\ex@
 \lineskip3\ex@\lineskiplimit3\ex@\halign\bgroup
 &\hfill$\m@th##$\hfill\cr}
\def\@endCD{\cr\egroup\egroup}
% Horizontal arrows with "sliding" length
\def\>#1>#2>{\amp@rs@nd@\setbox\z@\hbox{$\scriptstyle
 \;{#1}\;\;$}\setbox\@ne\hbox{$\scriptstyle\;{#2}\;\;$}\setbox\tw@
 \hbox{$#2$}\ifCD@
 \global\bigaw@\minCDaw@\else\global\bigaw@\minaw@\fi
 \ifdim\wd\z@>\bigaw@\global\bigaw@\wd\z@\fi
 \ifdim\wd\@ne>\bigaw@\global\bigaw@\wd\@ne\fi
 \ifCD@\hskip.5em\fi
 \ifdim\wd\tw@>\z@
 \mathrel{\mathop{\hbox to\bigaw@{\rightarrowfill}}\limits^{#1}_{#2}}\else
 \mathrel{\mathop{\hbox to\bigaw@{\rightarrowfill}}\limits^{#1}}\fi
 \ifCD@\hskip.5em\fi\amp@rs@nd@}
\def\<#1<#2<{\amp@rs@nd@\setbox\z@\hbox{$\scriptstyle
 \;\;{#1}\;$}\setbox\@ne\hbox{$\scriptstyle\;\;{#2}\;$}\setbox\tw@
 \hbox{$#2$}\ifCD@
 \global\bigaw@\minCDaw@\else\global\bigaw@\minaw@\fi
 \ifdim\wd\z@>\bigaw@\global\bigaw@\wd\z@\fi
 \ifdim\wd\@ne>\bigaw@\global\bigaw@\wd\@ne\fi
 \ifCD@\hskip.5em\fi
 \ifdim\wd\tw@>\z@
 \mathrel{\mathop{\hbox to\bigaw@{\leftarrowfill}}\limits^{#1}_{#2}}\else
 \mathrel{\mathop{\hbox to\bigaw@{\leftarrowfill}}\limits^{#1}}\fi
 \ifCD@\hskip.5em\fi\amp@rs@nd@}
% Rectangular commutative diagrams with diagonal arows
\newenvironment{CDS}{\@CDS}{\@endCDS}
\def\@CDS{\def\A##1A##2A{\llap{$\vcenter{\hbox
 {$\scriptstyle##1$}}$}\Big\uparrow\rlap{$\vcenter{\hbox{%
$\scriptstyle##2$}}$}&}%
\def\V##1V##2V{\llap{$\vcenter{\hbox
 {$\scriptstyle##1$}}$}\Big\downarrow\rlap{$\vcenter{\hbox{%
$\scriptstyle##2$}}$}&}%
\def\={&\hskip.5em\mathrel
 {\vbox{\hrule width\minCDaw@\vskip3\ex@\hrule width
 \minCDaw@}}\hskip.5em&}
\def\verteq{\Big\Vert&}
\def\novarr{&}
\def\noharr{&&}
\def\SE##1E##2E{\slantedarrow(0,18)(4,-3){##1}{##2}&}
\def\SW##1W##2W{\slantedarrow(24,18)(-4,-3){##1}{##2}&}
\def\NE##1E##2E{\slantedarrow(0,0)(4,3){##1}{##2}&}
\def\NW##1W##2W{\slantedarrow(24,0)(-4,3){##1}{##2}&}
\def\slantedarrow(##1)(##2)##3##4{%
\thinlines\unitlength1pt\lower 6.5pt\hbox{\begin{picture}(24,18)%
\put(##1){\vector(##2){24}}%
\put(0,8){$\scriptstyle##3$}%
\put(20,8){$\scriptstyle##4$}%
\end{picture}}}
\def\vspace##1{\noalign{\vskip##1\relax}}\relax\iffalse{%
\fi\let\amp@rs@nd@&\iffalse}\fi
 \CD@true\vcenter\bgroup\relax\iffalse{%
\fi\let\\=\cr\iffalse}\fi\tabskip\z@skip\baselineskip20\ex@
 \lineskip3\ex@\lineskiplimit3\ex@\halign\bgroup
 &\hfill$\m@th##$\hfill\cr}
\def\@endCDS{\cr\egroup\egroup}
% Triangular commutative diagrams
\newdimen\TriCDarrw@
\newif\ifTriV@
\newenvironment{TriCDV}{\@TriCDV}{\@endTriCD}
\newenvironment{TriCDA}{\@TriCDA}{\@endTriCD}
\def\@TriCDV{\TriV@true\def\TriCDpos@{6}\@TriCD}
\def\@TriCDA{\TriV@false\def\TriCDpos@{10}\@TriCD}
\def\@TriCD#1#2#3#4#5#6{%
\setbox0\hbox{$\ifTriV@#6\else#1\fi$}
\TriCDarrw@=\wd0 \advance\TriCDarrw@ 24pt
\advance\TriCDarrw@ -1em
\def\SE##1E##2E{\slantedarrow(0,18)(2,-3){##1}{##2}&}
\def\SW##1W##2W{\slantedarrow(12,18)(-2,-3){##1}{##2}&}
\def\NE##1E##2E{\slantedarrow(0,0)(2,3){##1}{##2}&}
\def\NW##1W##2W{\slantedarrow(12,0)(-2,3){##1}{##2}&}
\def\slantedarrow(##1)(##2)##3##4{\thinlines\unitlength1pt
\lower 6.5pt\hbox{\begin{picture}(12,18)%
\put(##1){\vector(##2){12}}%
\put(-4,\TriCDpos@){$\scriptstyle##3$}%
\put(12,\TriCDpos@){$\scriptstyle##4$}%
\end{picture}}}
\def\={\mathrel {\vbox{\hrule
   width\TriCDarrw@\vskip3\ex@\hrule width
   \TriCDarrw@}}}
\def\>##1>>{\setbox\z@\hbox{$\scriptstyle
 \;{##1}\;\;$}\global\bigaw@\TriCDarrw@
 \ifdim\wd\z@>\bigaw@\global\bigaw@\wd\z@\fi
 \hskip.5em
 \mathrel{\mathop{\hbox to \TriCDarrw@
{\rightarrowfill}}\limits^{##1}}
 \hskip.5em}
\def\<##1<<{\setbox\z@\hbox{$\scriptstyle
 \;{##1}\;\;$}\global\bigaw@\TriCDarrw@
 \ifdim\wd\z@>\bigaw@\global\bigaw@\wd\z@\fi
 \mathrel{\mathop{\hbox to\bigaw@{\leftarrowfill}}\limits^{##1}}
 }
 \CD@true\vcenter\bgroup\relax\iffalse{\fi\let\\=\cr\iffalse}\fi
 \tabskip\z@skip\baselineskip20\ex@
 \lineskip3\ex@\lineskiplimit3\ex@
 \ifTriV@
 \halign\bgroup
 &\hfill$\m@th##$\hfill\cr
#1&\multispan3\hfill$#2$\hfill&#3\\
&#4&#5\\
&&#6\cr\egroup%
\else
 \halign\bgroup
 &\hfill$\m@th##$\hfill\cr
&&#1\\%
&#2&#3\\
#4&\multispan3\hfill$#5$\hfill&#6\cr\egroup
\fi}
\def\@endTriCD{\egroup}
%%%%%%%%%%%%%%%  End of diagram macros.  %%%%%%%%%%%%%%%%%%%%%%%%%
% Skriptbuchstaben
\newcommand{\sA}{{\mathcal A}}
\newcommand{\sB}{{\mathcal B}}
\newcommand{\sC}{{\mathcal C}}
\newcommand{\sD}{{\mathcal D}}
\newcommand{\sE}{{\mathcal E}}
\newcommand{\sF}{{\mathcal F}}
\newcommand{\sG}{{\mathcal G}}
\newcommand{\sH}{{\mathcal H}}
\newcommand{\sI}{{\mathcal I}}
\newcommand{\sJ}{{\mathcal J}}
\newcommand{\sK}{{\mathcal K}}
\newcommand{\sL}{{\mathcal L}}
\newcommand{\sM}{{\mathcal M}}
\newcommand{\sN}{{\mathcal N}}
\newcommand{\sO}{{\mathcal O}}
\newcommand{\sP}{{\mathcal P}}
\newcommand{\sQ}{{\mathcal Q}}
\newcommand{\sR}{{\mathcal R}}
\newcommand{\sS}{{\mathcal S}}
\newcommand{\sT}{{\mathcal T}}
\newcommand{\sU}{{\mathcal U}}
\newcommand{\sV}{{\mathcal V}}
\newcommand{\sW}{{\mathcal W}}
\newcommand{\sX}{{\mathcal X}}
\newcommand{\sY}{{\mathcal Y}}
\newcommand{\sZ}{{\mathcal Z}}
% Sonderbuchstaben mit Doppellinie
\newcommand{\A}{{\mathbb A}}
\newcommand{\B}{{\mathbb B}}
\newcommand{\C}{{\mathbb C}}
\newcommand{\D}{{\mathbb D}}
\newcommand{\E}{{\mathbb E}}
\newcommand{\F}{{\mathbb F}}
\newcommand{\G}{{\mathbb G}}
\newcommand{\HH}{{\mathbb H}}
\newcommand{\I}{{\mathbb I}}
\newcommand{\J}{{\mathbb J}}
\newcommand{\M}{{\mathbb M}}
\newcommand{\N}{{\mathbb N}}
\renewcommand{\P}{{\mathbb P}}
\newcommand{\Q}{{\mathbb Q}}
\newcommand{\R}{{\mathbb R}}
\newcommand{\T}{{\mathbb T}}
\newcommand{\U}{{\mathbb U}}
\newcommand{\V}{{\mathbb V}}
\newcommand{\W}{{\mathbb W}}
\newcommand{\X}{{\mathbb X}}
\newcommand{\Y}{{\mathbb Y}}
\newcommand{\Z}{{\mathbb Z}}
\newcommand{\id}{{\rm id}}
\newcommand{\rank}{{\rm rank}}
%%%%%%%%%%%%%%%%%%%%%%%%%%%%%%%%%%%%%%%%%%%%%%%%%%%%%%%%%%%%%%
\title[Families of $K3$ surfaces]{Families of $K3$ surfaces
 over curves satisfying the  equality of Arakelov-Yau's type and Modularity}
%\title[]{} %mit Kurztitel in []
\author[Sun Xiaotao]{Xiaotao Sun}
\address{Department of Mathematics, The University of Hong Kong, 
Pokfulam Road, Hong Kong, and Institute of Mathemaics, 
Chinese Academy of Sciences, Beijing
100080, P. R. of China}
\email{xsun@math08.math.ac.cn}
\thanks{This work 
is supported by a  direct grant for Research
from the Chinese University of Hong Kong
(Project Code:  2060197).
Tan is also supported by
the 973 Foundation, the Foundation of EMC for Key Teachers and the PADS.
Sun is supported by a grant of NSFC for outstanding young researcher
(Project Code: 10025103).}
\author[Tan Sheng-Li]{Sheng-Li Tan}
\address{Department of Mathematics, East China Normal university,
Shanghai 200062, P. R. of China}
\email{sltan@euler.math.ecnu.edu.cn}
\author[Zuo Kang] {Kang Zuo}
\address{ Department of Mathematics, The Chinese University of Hong Kong,
Shatin, Hong Kong}
\email{kzuo@math.cuhk.edu.hk}
\maketitle
%%%%%%%%%%% Introduction %%%%%%%%%%%%%%%%%%%%%%%%%%%%%%%%%%%

Let $C$ denote  a smooth projective curve of genus $q$ over $\C,$ 
and $S'\subset C$ 
a finite set of points, 
and $f:  X^0\to C\setminus S'$ a smooth family of 
algebraic K3 surfaces, which extends to a family $f:  X\to C$ 
with semi-stable singular fibres over $S'.$ Let $S\subset S'$ 
denote the subset 
where the local monodromies of $R^2f_*\Bbb Z_{X^0}$ have infinite 
orders.
Let $\omega_{X/C}$ denote the dualizing sheaf. It is known that 
$f_*\omega_{X/C}$ is ample on 
$C $ by  Fujita if $f$ is not isotrivial \cite{F} 
 (See \cite{K} \cite{V} when the base of higher dimension). 
Jost and Zuo showed the following inequality \cite{J-Z}: 
\begin{equation}\label{arakelov_ineq_I}
\deg f_* \omega_{X/C}\leq \deg \Omega^1_C(\log S).
\end{equation}

If the iterated Kodaira-Spencer map of this  family is zero,  
one shows then a stronger 
inequality 
\begin{equation}\label{arakelov_ineq_II}
\deg f_*\omega_{X/C}\leq \frac{1}{2}\deg \Omega^1_C(\log S).
\end{equation}

These inequalities generalize the original Arakelov inequality 
$$\deg f_*\omega_{ X/C}\leq \frac{g}{2}\deg \Omega^1_C(\log S')$$
for a family of semi-stable curves of genus $g$.
The strict Arakelov inequality was proved in \cite {Tan} when 
$g\ge 2$ and $S'$ non-empty. 
If $S'$ is empty, Miyaoka-Yau inequality
implies a stronger inequality
$$\deg f_*\omega_{ X/C}\leq \frac{g-1}{6}\deg \Omega^1_C(\log S').$$
Thus the original
Arakelov inequality is always strict for $g\ge 2$. But there are 
families of Jacobians reaching the equality.
In general, Yau \cite{Y} proved the so-called Yau's Schwarz 
type inequality, 
which can be formulated as follows.
Let $(M, ds)$ be a Hermitian manifold with holomorphic sectional 
curvature bounded
above by a negative constant $K$, and let $(C\setminus S, ds_{\mu})$ 
be a Poincare
type metric. Then there exists a positive constant $c,$ such that 
for any holomorphic
map $\phi: C\setminus S\to M,$  one has 
$ \phi^*ds\leq c ds_{\mu}.$
It is the reason why we call the inequalities (\ref{arakelov_ineq_I}) 
and (\ref{arakelov_ineq_II}) 
are of Arakelov-Yau's type.
One  can choose $c=1$ in  (\ref{arakelov_ineq_I})  and  $c=1/2$ in 
(\ref{arakelov_ineq_II}). 
We will show in this note that both of them are optimal, 
and have a modularity meaning.\\

For a family  $f: E\to \P^1$ of non-constant semi-stable 
elliptic curves, Beauville proved
that $f$ has at least 4 singular fibres, which is 
equivalent to Arakelov inequality 
in this case. He obtained a complete classification 
for the families of semistable
elliptic curves when Arakelov inequality becomes equality. 
There are exactly $6$ such family
and $\Bbb P^1\setminus S'$ is a modular curve \cite{B2}.
In this note we shall study non-isotrivial algebraic families
of semi-stable K3 surfaces over curves when the inequality 
(\ref{arakelov_ineq_I}), 
or
(\ref{arakelov_ineq_II}) becomes an equality. 
The corresponding question has been considered in
\cite{V-Z} for families of abelian varieties. The final 
presentation of 
this note has been
influenced by \cite{V-Z}. It has also been motivated by
Mok's work on rigidity theorems of locally 
Hermitian symmetric spaces \cite{Mo1} and \cite{Mo2}, 
where he use the
Gaussian curvature of the induced metric on a 
holomorphic curves in a 
locally  
Hermitian symmetric space to characterize when this 
curve will be a 
totally geodesic embedding. \\

To state the main result, we recall some notation. 
Let $a:A^0\to C^0$ be a family of abelian surfaces 
with a section, then the
desingularization $Z^0\to C^0$ of the quotient 
$A^0/\{\pm 1\}\to C^0$ is a family
of Kummer surfaces (the so called Kummer construction). 
The rational map 
$A^0\to Z^0$ is called a rational
quotient of $A^0$. The family $a:A^0\to C^0$ is called 
the associated family
of abelian surfaces of $Z^0\to C^0$. In general, it is not 
true that every
family of Kummer surfaces has an associated family of 
abelian surfaces.
An involution $\imath$ on a $K3$ surface $X$ is called a 
Nikulin involution if
$\imath^*\omega=\omega$ for every 
$\omega\in H^0(X,\Omega^2_X)$. It is known
(Nikulin  \cite{N}) that every Nikulin involution 
$\imath$ has 
eight isolated fixed
points, and the rational quotient $X\to Z$ by 
$\imath$ is a $K3$ surface.

\begin {theorem} Let  $f: X\to C$  be a 
family of semi-stable $K3$ 
surfaces over $C,$
and suppose that $R^2f_*(\Bbb Z_{X^0})$ 
has infinite local 
monodromies at a non-empty set 
$S\subset C.$
If this family reaches the Arakelov bound in   
(\ref{arakelov_ineq_I}). Then we have

{\bf a)}  The general fibres of $f: X\to C$ have 
Picard number at least 19.

{\bf b)} After passing through a finite 
{\'e}tale cover $\sigma: C'\to C$, 
there exist an open set $C^0\subset C'$ and a global
Nikulin involution $\imath$ on 
$f:X^0=f^{-1}(C^0)\to C^0$ such that  
the rational quotient 
$X^0\to Z^0$ by $\imath$ is a family of 
Kummer surfaces over $C^0$, which has
an associated family of abelian surfaces that is
isogenous to the square product of a family of 
elliptic curves $g:E\to C'$. 
The projective  monodromy 
representation of the local system $R^1g_*(\Z_{E^0})$ 
extends to 
$$\tau: \pi_1(C'\setminus\sigma^{-1}S,*)\to PSL_2(\Z)$$ 
such that
$$C'\setminus\sigma^{-1}S\simeq \sH/\tau 
\pi_1(C'\setminus\tau^{-1}S,*).$$
\end{theorem}

A family of $K3$ surfaces satisfying Property b) 
will be called coming from
Nikulin-Kummer construction of the square 
product of a family of elliptic curves.

\begin{theorem}
If this family has the zero iterated Kodaira-Spencer map, 
and reaches
the Arakelov bound in (\ref{arakelov_ineq_II}).
Then the general fibres of $f: X\to C$ have the 
Picard number at least 18, 
after passing through a finite {\' e}tale cover 
$\sigma: C'\to C$, the monodromy
representation $\rho$ of $R^2f_*(\Z_{X^0})$ 
is of the form
$$\rho=\makebox{rank-2 trvial representation}\otimes
(\tau:\pi_1(C'\setminus\sigma^{-1}S,*)\to SL_2(\Z)),$$ 
$$and\qquad C'\setminus\sigma^{-1}S\simeq \sH/
\tau\pi_1(C'\setminus\tau^{-1}S,*).$$
\end{theorem}

\begin{remark}i) Theorem 0.1 can be used to explain 
the observation of 
B. Lian and S.-T. Yau (\cite{Li-Y1}, \cite {Li-Y2}) 
that the weight-2 VHS 
attached to $f$  of certain one dimensional 
families of $K3$ surfaces coming 
from the Mirror of $K3$ surfaces of Picard 
number $\geq 1$ can be expressed 
as the square products of the weight 1 VHS 
attached to some modular families 
of elliptic curves (also see \cite{Do}). 
Note that such a family must reach the
Arakelov bound in (\ref{arakelov_ineq_I}). 
We thank A. Todorov for  pointing 
out that to us. Note that, if $S=\emptyset$ 
then there is  another type 
families of $K3$ surfaces reaching the Arakelov 
bound (\ref{arakelov_ineq_I}). 
Namely, let $a:A\to C$ be a modular family of 
false elliptic curves, 
i.e. abelian surface whose endomorphism ring is 
isomorphic to an order 
of an indefinite quaternion algebra over $\Q$ 
(\cite{Sha}). Then the Kummer 
construction gives rise to a family $f:X\to C$ 
of smooth $K3$ surfaces reaching 
the Arakelov bound (\ref{arakelov_ineq_I}), 
and $C$ is a Shimura curve. 
One likes to know what is the mirror 
pair of this family. \\

ii) For a family $f: X\to C$ in Theorem 0.2 
one can find a family
$f': X'\to C,$ which comes from the 
Nikulin-Kummer construction of a product of 
a modular family of elliptic curves $g: E_1\to C$ 
with an elliptic curve $E_2$ 
over $\Bbb C, $ and 
such that sub VHSs of transcendental 
lattices of $f$ and $f'$ are
Hodge isometric to each other. 
Are there more closer geometric relations among
these families ?\end{remark}

Let $f:X\to \Bbb P^1$ be a Calabi-Yau 3-fold fibred 
by non-constant semi-stable
$K3$ surfaces. The triviality of $\omega_X$ implies that 
$\deg f_*\omega_{X/\Bbb P^1}=2.$

\begin{corollary} Let $f:X\to \Bbb P^1$ be a 
Calabi-Yau 3-fold fibred
by non-constant semi-stable $K3$ surfaces. 
Then the followings hold true:

{\bf i)} If the iterated Kodaira-Spencer 
map of $f$ is non-zero, 
then $f$ has at least 4 singular fibres. 
If $f$ has 4 singular fibres, 
then $X$ is rigid and birational to the 
Nikulin-Kummer
construction of a square product of a family 
of elliptic curves  
$g:E\to \Bbb P^1$. After passing through 
(if necessary) a double cover 
$E'\to E$, the family $g':E'\to \Bbb P^1$ 
is a modular family
of elliptic curves from the Beauville's 6 examples.

{\bf ii)} If the iterated Kodaira-Spencer 
map of $f$ is zero, 
then $f$ has at least 6 singular fibres. 
If $f$ has 6 singular fibres over 
$S\subset \Bbb P^1$, then $X$ is non-rigid, the  
general fibres have Picard number at least 18, and
$\Bbb P^1\setminus S\simeq \sH/\Gamma,$
where $\Gamma$ is a subgroup of $SL_2(\Z)$ of index 24.

\end{corollary}

\begin{remark} i) Any $K3$-fibred Calabi-Yau 
3-fold $f: X\to \Bbb P^1$ in i) is rigid
because of the modular construction for $X.$ 
Since all 6 examples of Beauville
are defined over $\Z$, we may assume that $X$ 
has a suitable integral model.
The $L-$series of $X$ is defined to be 
the $L-$series of the Galois representation on 
$H^3_{et}(\bar X,\Q).$ One should be able 
verify the so-called modularity conjecture
for $X$ as M-H Saito and N. Yui checked 
one example in \cite{S-Y}. That is, 
up to a finite Euler factor, 
$L(X,s)=L(f,s)$ for $ f\in S_4(\Gamma_0(N)).$

ii) Does any rigid Calabi-Yau 3-fold 
fibred by semi-stable $K3$ surfaces come from
the modular construction in i)?

iii) One can construct an example for the case ii) 
of Corollary 04.
Let $$g: E(4)\to X(4)$$ be
the modular family of elliptic curves corresponding 
to the congruence group $\Gamma(4).$
Then $X(4)\simeq \Bbb P^1$ with six cusps, and 
$\deg g_*\omega_{E(4)/X(4)}=2.$  
The Nikulin-Kummer construction to the product 
of $g: E(4)\to X(4)$ with a constant 
family of elliptic curves
gives an $K3$ fibred Calabi-Yau 3-fold 
reaching the upper bound in 
(\ref{arakelov_ineq_II}), which is non rigid. 

\end{remark}

{\bf Acknowledgment}
We owe the paper to discussions with H. Esnault,
B. Hassett, N. Mok, A. Todorov, E. Viehweg, S-T. Yau. 
Also, the conversation with S-W. Zhang
on Shimura curves is very helpful for us. Viehweg read through 
the  preliminary version, and made many valuable suggestions to 
improvment of this paper. We would like to thank all of them.

\section{Weight-2 VHS and $\R-$Splitting}

 Let $f: X\to C$ be a family of semi-stable K3 surfaces. 
Consider its weight-2 variation of Hodge structure 
(VHS for simplicity)
$$ \V=R^2f_*(\Z_{X^0}).$$

Let $S\subset C$ denote the subset, where the 
local monodromies of
$\Bbb V$ have infinite order. One has the canonical 
extension of Hodge
bundles 
$$ E^{p,q}=R^qf_*(\Omega^p_{ X/C}(\log\Delta),
\quad p+q=2,$$ 
together with the cup product of Kodaira-Spencer map
$$\theta^{p,q}: E^{p,q}\to E^{p-1,q+1}\otimes
\Omega^1_C(\log S).$$
$\theta=\theta^{2,0}+\theta^{1,1}$ is called t
he Higgs field of $\V$.

\begin{lemma} We have $\deg E^{2,0}\leq 
\deg\Omega^1_C(\log S)$, and if the equality 
$$\deg E^{2,0}=\deg\Omega^1_C(\log S)$$ 
holds, then there is a  real splitting
$ \V\otimes\R=\W\oplus\U,$
which is orthogonal w.r.t. the polarization, 
and $ \U$ is unitary.
The corresponding Higgs bundle splitting is 
$$ (E^{2,0}\oplus E^{1,1}_1\oplus E^{0,2},\theta)
\oplus (E^{1,1}_2,0)$$
where $E^{1,1}=E^{1,1}_1\oplus E^{1,1}_2$
and $ E^{1,1}_1$ is a line bundle of 
degree zero such that 
$$ \theta: E^{2,0}\to E_1^{1,1}\otimes 
\Omega^1_C(\log S),\quad
\theta: E_1^{1,1}\to E^{0,2}\otimes 
\Omega^1_C(\log S)$$
are isomorphisms.
\end{lemma}

\begin{proof}   Consider the map
$\theta^{1,1}: E^{1,1}\to E^{0,2}\otimes
\Omega^1_C(\log S),$
let $E_2^{1,1}\subset E^{1,1}$ denote the 
kernel of $\theta^{1,1},$
then $(E_2^{1,1},0)$ is a Higgs sub-bundle.\\

{\bf Claim:}\quad  $\deg E^{1,1}_2\leq 0,$ and if
the equality holds then the Higgs subbundle
$$ (E^{1,1}_2, 0)\subset (E,\theta) $$
induces a splitting
$(E,\theta)=(E^{2,0}\oplus E^{1,1}_1\oplus E^{0,2},
\theta)\oplus (E^{1,1}_2, 0),$
which corresponds to a splitting of the local 
system over $\C$
$ \V\otimes \C=\W\oplus \U.$   \\

{\bf Proof of the claim:}\quad  Let $h$ denote 
the Hodge metric on $E|_{C\setminus S},$
and let $\Theta(E|_{C\setminus S},h)$ be 
its curvature form. Then we have
(\cite{G}, Chapter II) 

$$\Theta(E|_{C\setminus S})+
\theta\wedge\bar\theta+\bar\theta\wedge\theta=0, $$
where $\bar \theta$ is the complex conjugation 
of $\theta$ with respect to
$h.$ Consider the 
$ \mathcal C^{\infty}$-orthogonal (for $h$)
decomposition 
$ E|_{C\setminus S}=E^{1,1}_2|_{C\setminus S}\oplus 
{E^{1,1}_2|_{C\setminus S}}^{\perp},$ one has  
$$ \Theta(E^{1,1}_2|_{C\setminus S},h)=
\Theta(E|_{C\setminus S},h)|_{E^{1,1}_2}+
\bar A\wedge A=
-(\theta\wedge\bar \theta)|_{E^{1,1}_2}-
(\bar\theta\wedge\theta)|_{E^{1,1}_2}+
\bar A\wedge A,$$
where $ A\in A^{1,0}({Hom}(E^{1,1}_2, {E^{1,1}_2}^{\perp}))$ 
is the second fundamental
form of the subbundle $ E^{1,1}_2\subset E,\,$ and $\, \bar A$ 
is the complex conjugation
with respect to $ h.$ 
Since $ \theta (E^{1,1}_2)=0,$ we have 
$(\bar\theta\wedge\theta)|_{E^{1,1}_2}=0.$ Hence
$$ \Theta(E^{1,1}_2|_{C\setminus S'},h)
=-(\theta\wedge\bar\theta)_{E^{1,1}_2}+\bar A\wedge A.$$
$ \Theta(E^{1,1}_2|_{C\setminus s'} ,h)$ is 
negative semidefinite
since $ \theta\wedge\bar\theta_{E^{1,1}_2} $ 
is positive semidefinite 
and $ \bar A\wedge A$ is negative semidefinite. 
Since the local monodromies around points in $S$ 
are unipotent,   
 $ {Tr}\,\Theta(E^{1,1}_2|_{C\setminus S'} ,h)$ 
represents (by \cite{Sch}) the Chern class 
$c_1(E^{1,1}_2)$ as a current. Thus
$$\deg E^{1,1}_2=\int_{C\setminus S}{Tr}
\Theta(E^{1,1}_2|_{C\setminus S} ,h)\leq 0,$$
and $\Theta(E^{1,1}_2|_{C\setminus S},h)=0$ 
if $\deg E^{1,1}_2=0.$ 
This implies that $\bar \theta(E^{1,1}_2)=0$ 
and $ A=0.$
Altogether show that the sub-Higgs bundle 
$(E^{1,1}_2,0)$ of $(E,\theta)$ induces a   
splitting of the Higgs bundle
$$ (E,\theta)=(E^{2,0}\oplus E^{1,1}_1
\oplus E^{0,2},\theta)\oplus (E^{1,1}_2,0)$$
and the corresponding splitting 
$\V\otimes\C=\W\oplus \U$ of the complex 
local system. Thus the claim is proved.
\\

Let $I\subset E^{0,2}\otimes
\Omega^1_C(\log S)$ be the image of
$\theta^{1,1},$ then, by the exact sequence
$$ 0\to E^{1,1}_2\to E^{1,1}\to I\to 0 ,$$
and note that $\deg E^{1,1}=0,$ one gets
$$ -\deg E^{2,0}+\deg\Omega^1_C(\log S)=
\deg(E^{0,2}\otimes\Omega^1_C(\log S))
\geq \deg I=-\deg E^{1,1}_2\geq 0.$$
Hence $\deg E^{2,0}\leq \deg\Omega^1_C(\log S)$ 
and the equality holds if and only
if $\deg E^{1,1}_2=0$ and $I=E^{0,2}
\otimes\Omega^1_C(\log S),$ which is our $E^{1,1}_1.$
It is easy to see that the Higgs field of $\W$ 
is an isomorphism, thus $\W$ is irreducible
over $\C$. 
Now we only need to show that the 
decomposition $\V\otimes\C=\W\oplus \U$
can be, in fact, defined over $\R.$ 
Taking the complex conjugation on $\W$
one has
$$\overline{\W}\subset \overline{\V\otimes \C}= 
\V\otimes \C.$$

$\bar\W$ is again of the Hodge 
type $(2,0)+(1,1)+(0,2)$, irreducible
and with the non-zero Higgs field.
The projection 
$ p: \overline{\W}\subset  \V\otimes \C \to \U$
can not be injective since $\U$ is unitary. 
Moreover, 
since $\bar \W$ can not have a proper sub local system,
this projection must be zero. Thus $\overline{\W}=\W$
and we obtain a real sub local system
$ \W\subset \V\otimes \R.$
The intersection form restricted to  
$\W$ is non-degenerated.
Thus the orthogonal complement of $\W$ 
with respect to the 
intersection form gives the desired 
real decomposition
$ \V\otimes \R=\W\oplus \U.$ 
\end{proof}

\begin{lemma}  If the iterated Kodaira-Spencer map 
$\theta^{1,1}\theta^{2,0}=0,$ then
$$\deg E^{2,0}\leq \frac{1}{2}\deg\Omega^1_C(\log S).$$ 
When the equality
$\deg E^{2,0}=\frac{1}{2}\deg\Omega^1_C(\log S)$
holds, then there is a real splitting
$$\V\otimes\R=\W\oplus\U,$$
which is orthogonal w.r.t. the polarization, 
and $\U$ is unitary.
The corresponding Higgs bundle splitting is
$$(E^{2,0}\oplus (E^{1,1}_1\oplus {E^{1,1}_1}^*)
\oplus E^{0,2},\theta)\oplus (E_2^{1,1},0)$$
where $E^{1,1}_1$ and ${E^{1,1}_1}^*$  are  
sub line bundles of  $E^{1,1}$ with 
$$\deg E^{1,1}_1=-\deg E^{2,0}=
-\frac{1}{2}\deg\Omega^1_C(\log S),$$
and $ E^{1,1}=E^{1,1}_1\oplus 
{E^{1,1}_1}^*\oplus E^{1,1}_2.$ The Higgs field
$$ \theta: (E^{2,0}\oplus (E^{1,1}_1
\oplus {E^{1,1}_1}^*)\oplus E^{0,2})\to  
(E^{2,0}\oplus (E^{1,1}_1\oplus {E^{1,1}_1}^*)
\oplus E^{0,2})\otimes\Omega^1_C(\log S)$$
is defined by $\theta=\tau\oplus-\tau^*,$ where
$\tau: E^{2,0}\simeq E^{1,1}_1
\otimes\Omega^1_C(\log S),\quad
E^{1,1}_1\to 0.$
\end{lemma}

\begin{proof}\quad since $\theta^{1,1}\theta^{2,0}=0,$ 
the map $\theta^{2,0}$ factors through
$$ \theta^{2,0}: E^{2,0}\to E^{1,1}_1\otimes
\Omega^1_C(\log S),$$
where $E^{1,1}_1\subset E^{1,1}$ is a sub-line bundle 
such that $\theta^{1,1}(E^{1,1}_1)=0.$
Thus 
$$(E^{2,0}\oplus E^{1,1}_1,\theta^{2,0})
\subset (E,\theta)$$
is a rank-2 Higgs sub bundle. By the same arguments 
in the proof of Lemma 1.1, one has
$\deg E^{2,0}\oplus E^{1,1}_1\leq 0,$ 
thus
$$\deg E^{2,0}\leq\frac{1}{2}\Omega^1_C(\log S).$$
If the equality holds, then
$ \theta^{2,0}=:\tau: E^{2,0}\to E^{1,1}_1
\otimes\Omega^1_C(\log S)$
is an isomorphism with
$\deg E^{1,1}_1=-\deg E^{2,0}=-{\frac{1}{2}}\deg\Omega^1_C(\log S),$
and the Higgs sub bundle
$(E^{2,0}\oplus E^{1,1}_1,\theta^{2,0})\subset (E,\theta)$
gives rise to a complex sub local system $\W_1\subset\V\otimes\C.$
The dual $\bar \W_1\subset\V\otimes\C$ corresponds to
Higgs subbundle 
$$(E^{2,0}\oplus E^{1,1}_1)^*={E^{1,1}_1}^*\oplus E^{0,2}$$
together with the Higgs field
$-\tau^*:{E^{1,1}_1}^*\to E^{0,2}\otimes\Omega^1_C(\log S).$
The sub-local system  $\W:=\W_1\oplus\bar\W_1$
is real, and the intersection form restricted to
$\W$ is non-degenerated. Hence, the orthogonal
complement defines the desired decomposition.

\end{proof}

\section{Splitting over $\bar\Q$}

We start with a very simple observation. Suppose that  $\V$ is 
a local system defined over $\bar\Q.$  Fixing a positive integer $r,$ 
let $\mathcal G(r,\V)$ denote the set of
all rank-r sub-local systems of $\V$. Then $\mathcal G(r,\V)$ is a 
projective variety
defined over $\bar\Q.$ The following property is well known.

\begin{lemma} If $[W]\in \mathcal G(r,\V)$ is an isolated point, 
then $W$ is defined over $\bar\Q.$
\end{lemma}

\begin{lemma} The $\R$-splittings $\V\otimes\R=\W\oplus\U$
in  Lemma 1.1 and Lemma 1.2 can be defined over 
$\bar \Q.$\end{lemma} 

\begin{proof} By Lemma 2.1, one only needs to show that $\W$ is 
a rigid sub-local system
of $\V\otimes \C$. Suppose that there is a family of sub-local systems
$$\{\W_t\},\quad \W_0=\W.$$
By semi-continuity, the Higgs fields $\theta^{p,q}$ of $\W_t$ are again 
isomorphisms for
$t$ being sufficiently closed to $0$. Then the projection
$ \W_t\to \V\otimes\C\to\U $
must be zero, otherwise, $\W_t$ would contain a non-trivial unitary 
component, which contradicts that $\theta^{p,q}$ are isomorphisms. 
Hence $\W_t=\W.$\\

Similarly, we show that the sub-local system
$\W=\W_1\oplus\bar\W_1\subset\V=\W\oplus \U$
is rigid. Suppose that there is a family of sub local systems
$\{\W_t\}$ with $\W_0=\W$,
we decompose $\W_t$ into the direct sum of irreducible components over
$\C,$ which has only following possible types up to isomorphism
$$ \W_1\oplus\bar \W_1;\quad \W_1\oplus\U';\quad \bar\W_1\oplus\U'';\quad
\U''',$$
where $\U',\,\U'',\, \U'''$ are unitary. By semicontinuity,
the last three cases are
impossible if $t$ is sufficiently closed to $0$
(otherwise $\theta^{1,1}$ would be zero). Thus
$$ \W_t\simeq \W_1\oplus\bar \W_1,$$
which implies that the projection
$\W_t\to\V\otimes \C\to\U $ 
must be zero. Otherwise, $\W_1$ would contain a non-trivial unitary 
component, which contradicts that the Higgs fields of $\W$ are 
isomorphisms.
\end{proof}

\section{Splitting over $\Q$ and $\Z$-structures}

We call the splitting in Lemma 1.1 of type (\ref{arakelov_ineq_I}) and
the splitting in Lemma 1.2 of type (\ref{arakelov_ineq_II}).

\begin{lemma}If $S\not=\emptyset$, the splittings in Lemma 2.2 can be 
defined over $\Q.$\end{lemma}

\begin{proof} Let $\V\otimes K=\W\oplus\U $
be the splitting of type (\ref{arakelov_ineq_I}) in Lemma 2.2, 
where $K$ is a Galois extension of 
$\Q$. For any $\sigma\in Gal(K/\Q)$, we claim that $\sigma\W=W.$ 
Otherwise, the projection
$p: \sigma\W\to\V\otimes K\to \U$ 
must be non-zero and
$\sigma\W$ is isomorphic to a unitary sub local system $\U'\subset \U$ 
under $p$
since $\W$ is irreducible (thus $\sigma\W$ is also irreducible).
Let $\gamma$ be a short loop around $s\in S.$ Then the monodromy matrix
$\rho_{\W}(\gamma)$ has infinite order, hence $\rho_{\sigma\W}(\gamma)$ 
has also infinite order, which contradicts that $\rho_{\U'}(\gamma)$ is 
identity.
We proved that $\W$ is invariant under $Gal(K/\Q).$ 
Hence $\W$ is defined over $\Q$ and
the orthogonal complement of $\W\subset\V\otimes\Q$ w.r.t. 
the intersection form defines an $\Q-$splitting
$$\V\otimes\Q=\W\oplus\U.$$

By the same argument, we show that the splitting of 
type (\ref{arakelov_ineq_II}) in lemma 2.2
is also defined over $\Q.$\end{proof}

\begin{lemma} After passing through a finite etale cover of $C$ the 
splittings of
type (\ref{arakelov_ineq_I}) and (\ref{arakelov_ineq_II}) in Lemma 3.1 
induce $\Z-$sub lattices
$$ \V\supset \W_{\Z}\oplus\Z^{19},\quad \V\supset \W_{\Z}\oplus\Z^{18} $$
such that
$\V\otimes\Q=(\W_{\Z}\oplus\Z^{19})\otimes\Q$
and 
$\V\otimes\Q= (\W_{\Z}\oplus\Z^{18})\otimes\Q,$   
where $\Z^{19}$, $ \Z^{18}$ is respectively a rank-19 constant $\Z$-lattice 
of type-(1,1) 
and a rank-18 constant $\Z$-lattice of type-(1,1).\end{lemma}

\begin{proof} Let $\W_{\Z}=\V\cap\W,\quad \U_{\Z}=\V\cap\U.$
It is easy to check that
$$ \W_{\Z}\otimes\Q=\W,\quad \U_{\Z}\otimes\Q=\U, $$ 
thus $\W_{Z}$ and $\U_{\Z}$ are lattices in $\W$ and $\U$.
Since $\U$ is unitary and carries an $\Z-$structure, 
the monodromy group of $\U$ 
is finite. Since the local monodromies of $\U$ around $S$ are trivial, 
$\U$ extends to a local system on $C.$
Therefore, after passing through the cover corresponding to this monodromy 
group, $\U$ becomes a constant local system $\Z^{19}$, $\Z^{18}$ 
respectively.\end{proof}

\begin{corollary} Let $f:X\to C$ be a family of semi-stable $K3$ surfaces over a curve
$C.$ When it reaches the upper bound
$\deg f_*\omega_{X/C}=\deg\Omega^1_C(\log S),$
then the Picard number of the general fibres is at least $19$.
If $\theta^{1,1}\theta^{2,0}=0$ and $f$ reaches the upper bound 
$\deg f_*\omega_{X/C}=\frac{1}{2}\deg\Omega^1_C(\log S),$
then the Picard number of the general fibres is at least $18$.
\end{corollary}

\section{Nikulin and Kummer construction}

Let $f:X\to C$ be a family of semi-stable  $K3$ surfaces, which reaches
the upper  bound
$deg f_*\omega_{ X/C}=\deg\Omega^1_C(\log S).$
By Lemma 3.2, after passing through a finite {\'e}tale cover of $C,$ 
one has
$$ R^2f_*(\Z_{X^0})\otimes\Q=\W\oplus\Q^{19},$$
where $\W$ is an $\C$-irreducible representation of $\pi_1(C\setminus S,*)$ 
and $\Q^{19}$ is a constant local system of rank 19 such that
$ \Q^{19}_t\subset NS(X_t)\otimes\Q$ for any $t\in C\setminus S.$
We obtain therefore,

\begin{lemma} For any $t\in C\setminus S$,  the Picard number
$\rho(X_t)\geq 19$
and for any class $s_t\in \Q^{19}_t\subset Pic(X_t)\otimes\Q$
there is a $\Q$-divisor
$D\in Div(X)\otimes\Q$ such that $D|_{X_t}=s_t$.\end{lemma}

Let $Y$ be an algebraic K3 surface and 
$H^2(Y,\ Z)=T_Y\oplus NS(Y)$
be the orthogonal decomposition. $T_Y$ is the so called transcendental
lattice of $Y$, which is even and  has signature $(2,\,20-\rho(Y))$.
It is well-known that as lattices
$$H^2(Y,\Z)\cong U^3\oplus E_8(-1)^2.$$
We recall some results about embeddings of lattices (see \cite{Mo}  
and references there)

\begin{lemma}(Theorem 2.4 of \cite{L-P}, or see Corollary 2.6 of [Mo])
Let $T$ be a non-degenerate even lattice of rank $r$. Then there is
a primitive embedding $$T\hookrightarrow U^r$$
In particular, if $\rho(X)\ge 19$, then
there is a primitive embedding $$T_X\hookrightarrow U^3.$$\end{lemma}

\begin{lemma} If $12<\rho\le 20$, then every even lattice $T$ of
signature $(2,20-\rho)$ occurs as the transcendental lattice of some
algebraic K3 surface and the primitive embedding
$T\hookrightarrow U^3\oplus E_8(-1)^2$
is unique.\end{lemma}

\begin{theorem}(\cite{Mo}) If $\rho(Y)\geq 19$, then there exists a
primitive embedding
$$\varphi:E_8(-1)^2\hookrightarrow NS(Y)\subset H^2(Y,\Z)$$
and a Nikulin involution $\tau: Y\to Y$ such that 
$\tau^*:H^2(Y,\Z) \to H^2(Y,\Bbb Z)$ is identity on 
$(\varphi(E_8(-1)^2)^{\bot}$.\end{theorem}

\begin{proof} By Lemma 4.2, there is a primitive embedding
$\phi:T_Y\hookrightarrow U^3,$ 
thus a primitive embedding
$\phi\oplus 0:T_Y\hookrightarrow U^3\oplus E_8(-1)^2.$
By Lemma 4.3 (uniqueness), the above embedding is isomorphic to
$$T_Y=NS(X)^{\bot}\subset H^2(Y,\Z)\cong U^3\oplus E_8(-1)^2.$$
Thus, there is a primitive embedding 
$$\psi:E_8(-1)^2\hookrightarrow
T_Y^{\bot}=NS(Y)\subset H^2(Y,\Z).$$
Let $\{c_j^1\}_{1\le j\le 8}$ and $\{c_j^2\}_{1\le j\le 8}$ be the
bases of $E_8(-1)\oplus 0$ and $0\oplus E_8(-1)$ and
$$g:H^2(Y,\Z)\to H^2(Y,\Z)$$
be defined as: 
$g(\psi(c^1_j))=\psi(c^2_j)$, $\,g(\psi(c^2_j))=\psi(c^1_j)$ and
$\, g(e)=e$ for any $e\in (\psi(E_8(-1)^2))^{\bot}$.  
Then, by theorems of Nikulin (see Theorem 5.6 of [Mo]), 
there is a Nikulin involution
$\tau:Y\to Y$ and $w\in W(Y)$ (the group of Picard-Lefschetz reflections)
such that $\tau^*=w\cdot g\cdot w^{-1}$. Let
$$\varphi:E_8(-1)^2\stackrel{\psi}{\to}H^2(Y,\Z)\stackrel{w}{\to}H^2(Y,\Z),$$
then 
$\varphi:E_8(-1)^2\hookrightarrow NS(Y)\subset H^2(Y,\Bbb Z)$ is
another primitive embedding, and
$$\tau^*(\varphi(c_j^1))=\varphi(c^2_j),\quad
\tau^*(\varphi(c_j^2))=\varphi(c^1_j),\quad \tau^*(e)=e,\quad
\forall e\in(\varphi(E_8(-1)^2))^{\bot}.$$
\end{proof}

Let $t_0\in C\setminus S$ be a point such that the fibre $X_{t_0}$ 
satisfying $\rho(X_{t_0})=19$. Thus,
$$\Q^{19}_{t_0}= NS(X_{t_0})\otimes\Q.$$
Since the monodromy action of $\pi_1(C\setminus S, t_0)$ on 
$\Q^{19}_{t_0}$ is trivial,
$\varphi(c^1_j)$ and $\varphi(c^2_j),\,1\leq j\leq 8$ can be lifted 
to divisors
 $D^1_j$
and $D^2_j,\, 1\leq j\leq 8$ on $X.$ Then we have

\begin{lemma} For any $t\in C\setminus S$, let ${d^i_j}_t=D^i_j|_{X_t}
\in H^2(X_t,\Z)$. Then $\{{d^i_j}_t\}_{1\le j\le 8}$ ($i=1,2$)
generate a sublattice of $H^2(X_t,\Z)$, which is isomorphic to
$E_8(-1)^2$ such that $E_8(-1)^2\hookrightarrow H^2(X_t,\Z)$ is a
primitive embedding, $E_8(-1)\oplus 0$ and $0\oplus E_8(-1)$ are isomorphic
to $\Bbb Z\{{d^1_j}_t,j=1,...,8\}$ and $\Bbb Z\{{d^2_j}_t,j=1,...,8\}$
\end{lemma}

\begin{proof} The proof is straightforward. For example, to prove
that $\{{d^1_j}_t\}_{1\le j\le 8}$ are $\Z$-linearly independent:
if $\sum n_j{d_j^1}_t=0$ in $H^2(X_t,\Z)$, we claim that
$\sum n_j\varphi(c_j^1)=0,$ which will imply the $\Z$-linearly
independence of $\{{d^1_j}_t\}_{1\le j\le 8}$.
The claim is clear, otherwise there is a $A\in NS(X_{t_0})$ such that
$(\sum n_j\varphi(c_j^1),\,A)\not= 0.$ Let $\tilde A$ be a lifting of $A$, 
then
$$\aligned (\sum n_j{d_j^1}_t,\,\tilde A|_{X_t})&=
(\sum n_jD^1_j|_{X_t},\,
\tilde A|_{X_t})=(\sum n_jD^1_j|_{X_{t_0}},\,
\tilde A|_{X_{t_0}})\\&=(\sum n_j\varphi(c_j^1),\,A)\neq 0.\endaligned$$
To see that the embedding $E_8(-1)^2\hookrightarrow H^2(X_t,\Z)$
is primitive, let $B\in H^2(X_t,\Z)$ be a class with
$mB\in \Z\{{d^i_j}_t,i=1,2, j=1,...,8\}$. Then $B$ is invariant under the
monodromy, and thus there is a lifting $\tilde B$ of $B$. Since
$\varphi:E_8(-1)^2\hookrightarrow H^2(X_{t_0},\Z)$ is primitive
and $m\tilde B|_{X_{t_0}}\in\varphi(E_8(-1)^2)$,
$\tilde B|_{X_{t_0}}=\sum n^i_j\varphi(c_j^i)$. Then
$$(m(\tilde B-\sum n^i_jD^i_j)|_{X_t},\,
m(\tilde B-\sum n^i_jD^i_j)|_{X_t})=0$$
and $(m(\tilde B-\sum n^i_jD^i_j)|_{X_t}, H|_{X_t})=0$, which implies that
$m(\tilde B-\sum n^i_jD^i_j)|_{X_t}=0$
by Hodge index theorem since
$m(\tilde B-\sum n^i_jD^i_j)|_{X_t}=mB-m\sum n^i_j{d^i_j}_t$ is an
algebraic class. Thus $(\tilde B-\sum n^i_jD^i_j)|_{X_t}=0$ since
$H^2(\X_t,\Z)$ is torsion free.

\end{proof}

Let $E=\bigoplus_{p+q=2}E^{p,q}$ be the canonical extension of the 
Hodge bundle associated to the local system $R^2f_*(\Z_{X^0})$, and 
$\mathcal E nd(E)\to C$ be the endomorphism bundle over $C$, 
which represents the
functor 
$$\mathcal E nd(E)^{\sharp}:\{\text{schemes over $C$}\}\to\{sets\}$$
where $\mathcal E nd(E)^{\sharp}(T)=
\{\text{bundle morphism $E_T\to\E_T$ over $T$}\}$.
For $t\in C\setminus S$, by Lemma 4.5, we can define an isometric 
involution
$$g_t:H^2(X_t,\Z)\to H^2(\X_t,\Z)$$
by 
$g_t({d^1_j}_t)={d^2_j}_t$, $\, g_t({d^2_j}_t)={d^1_j}_t$,
$\,g_t(e)=e$ for all
$e\in\Z\{{d^i_j}_t\}^{\bot}$ and $1\leq j\leq 8.$
It is easy to see that 
$g_t:H^2(X_t,\Bbb Z)\to H^2(X_t,\Bbb Z)$ is a morphism of 
$\pi_1(C\setminus S)$-modules. Thus, they give rise an involution
$$g:R^2f_*(\Bbb Z_{X^0})\to R^2f_*(\Z_{X^0})$$
of local system, which corresponds to a section  
$g\in H^0(C\setminus S, \mathcal E nd(E))$.

\begin{lemma} The section $g\in H^0(C\setminus S, \mathcal E nd(E))$ 
defined above can be extended to
a section in $ H^0(C, \mathcal E nd(E)),$ and thus $g$ is an algebraic 
section.\end{lemma}

\begin{proof} Recall that 
$R^2f_*(\Bbb Z_{X^0})\otimes\Q=\W\oplus \Q^{19}$
and the canonical extension of the Hodge bundle corresponding to 
$R^2f_*(\Bbb Z_{X^0})$ can be written into
$$ (E,\theta)=(E_{\W},\theta)\oplus(\mathcal O_C^{19},0),$$
where $(E_{\W},\theta)$ and $ (\mathcal O_C^{19},0) $  are the canonical 
extension of the Hodge bundles corresponding to
$\W$ and $\Q^{19}$ respectively. By the construction of $g$, it is 
identity on ${\W}$ (thus extended to $E_{\W}$), and is well-defined on 
the constant lattice $\Z^{19}$. Thus it is clear that $g$ can be extended
on $C$.

\end{proof}

\begin{lemma} Let $H$ be an ample divisor on $X$ and 
$g_t:H^2(X_t,\Z)\to H^2(X_t,\Z)$
be the Hodge isometry involutions defined above. 
Then there exists a non-empty Zariski open set 
$C^0\subset C\setminus S$ such that $g_t(H|_{X_t})$ is an ample divisor 
for any $t\in C^0$,
In particular, $g_t$ is an effective Hodge isometry for any $t\in C^0$.
\end{lemma}

\begin{proof} We may write
$H|_{X_{t_0}}=\sum n_j^1\varphi(c_j^1)\,+ \sum n_j^2\varphi(c_j^2)\,\,+\,e,$ 
where $e\in\varphi(E_8(-1)^2)^{\bot}.$ Let $E$ be a lifting of $e$ and
$$D=\sum_{j=1}^8 n_j^1D_j^1\,+\sum_{j=1}^8 n_j^2D^2_j\,\,+E,\quad
\tilde D=\sum_{j=1}^8 n_j^1D_j^2\,+\sum_{j=1}^8 n_j^2D^1_j\,\,+E.$$
Then, for any $t\in C\setminus S$, $H|_{X_t}=D|_{X_t}$ and
$g_t(D|_{X_t})=\tilde D|_{ X_t}.$ Thus $D$ is a relative ample
divisor on $f^{-1}(C\setminus S)$ and $\tilde D|_{X_{t_0}}$ is ample 
(here we have choosen $t_0$ such that $g_{t_0}$ is effective). 
Thus there exists
a Zariski open set $C^0\subset C\setminus S$ such that 
$\tilde D$ is relative
ample on $f^{-1}(C^0)$.

\end{proof}

\begin{lemma} The $g$ induces an involution
$\tau: f^{-1}(C^0)\to f^{-1}(C^0)$ over $C^0$ such that
$\tau_t :  X_t\to X_t$ (for $t\in C^0$)
are Nikulin involutions with $\tau^*_t=g_t$.
\end{lemma}  

\begin{proof} Let $\mathcal L=D+\tilde D$, where $D$ and $\tilde D$ 
are the divisors
defined in the proof of Lemma 4.7. Then we know that $\mathcal L$ 
is relative
ample on $f^{-1}(C^0)$ and $\mathcal L_t=\mathcal L|_{X_t}$ is invariant 
under the involution $g_t$.
Let $\pi:Aut^{\mathcal L}(f^{-1}(C^0)/C^0)\to C^0$ denote the 
automorphism group scheme,
which represents the functor
$$Aut^{\mathcal L}_{f^{-1}(C^0)/C^0 }(T)=
\left\{\aligned &\text{Isomorphisms $h:f^{-1}(C^0)\times_{C^0}T\to
f^{-1}(C^0)\times_{C^0}T$}\\&\text{over $T$ such that
$h^*(p^*_T\mathcal L)=p^*_T(\mathcal L)$}
\endaligned\right\}.$$
Thus there exists a universal automorphism
$$ \begin{array}{cccccccccc}
f^{-1}(C^0)\times_{C^0}Aut^{\mathcal L}(f^{-1}(C^0)/C^0)&\stackrel {h}{\to}&
f^{-1}(C^0)\times_{C^0}Aut^{\mathcal L}(f^{-1}(C^0)/C^0)\\
\tilde f\downarrow &   & \tilde f\downarrow\\   
Aut^{\mathcal L}(f^{-1}(C^0)/C^0)&=&Aut^{\mathcal L}(f^{-1}(C^0)/C^0) \\
\end{array}$$
and $h^*$ induces an endomorphism $\pi^*E\to\pi^*E$, which gives a
homomorphism
$$\begin{array}{cccccccccc}
Aut^{\mathcal L}(f^{-1}(C^0)/C^0)&\stackrel{\alpha}{\to}& \mathcal E nd(E) \\   
\pi\downarrow &  & \downarrow\\
C^0&=&C^0.\\
\end{array}$$
By Torelli theorem of K3 surfaces, $\alpha$ is injective. On the other
hand, the fibres of $\alpha$ are isomorphic to group schemes, which
are smooth. Thus $\alpha$ is an embedding. By Lemma 4.6 and Lemma 4.7,  
$g(C^0)$ is
algebraic and  contained in the image of $\alpha$, which gives a section of
$\pi:Aut^{\mathcal L}(f^{-1}(C^0)/C^0)\to C^0.$
That is an automorphism
$$ \begin{array}{cccccccccc}
f^{-1}(C^0)&\stackrel{\tau}{\to}& f^{-1}(C^0)\\
f\downarrow & &f\downarrow\\
C^0 &=&C^0\\
\end{array}$$ 
such that $\tau^*_t=g_t$ for any $t\in C^0$. Thus $\tau_t$ are Nikulin
involutions, i.e. $\tau_t^*\omega=\omega$ for any
$\omega\in H^{2,0}(X_t).$

\end{proof}

Since all fibres $X_t$ are algebraic K3 surfaces, the $\tau_t$ gives
rise  a Shioda-Inose structure on $X_t$ by theorems of Morrison
(see Theorem 6.3 of \cite{Mo}). Let $g: Z^0\to C^0$ be the
desingularization of $f^{-1}(C^0)/\tau\to C^0$. Then $g: Z^0\to C^0$
is a family of Kummer surfaces and there exist divisors $N_1,...,N_8$ on
$Z^0$ such that their restrictions $(N_1)_t,...,(N_8)_t$ on
$Z^0_t$ are the exceptional $(-2)$-curves of the double points
of $X_t/\tau_t$ (produced by the eight isolated fixed points of
$\tau_t$). By Lemma 3.2, we write
$ R^2f_*(\Z_{f^{-1}(C^0)})=\W\oplus\Z^{19}$.
Then we have (see Lemma 3.1 of \cite{Mo})
$$R^2g_*(\Z_{Z^0})\simeq (\W\oplus{\Z^{19}}^{\tau})(2)\oplus
\Z[N_1,...,N_8],$$
where $ {\Z^{19}}^{\tau} $ is the invariant sub local system of $\Z^{19}$
under $\tau,$ $(\W\oplus{\Z^{19}}^{\tau})(2)$ has the same underlying local system as
$(\W\oplus{\Z^{19}}^{\tau}),$ and with an intersection  by multiplication
by 2 of the the intersection form on $(\W\oplus{\Z^{19}}^{\tau}).$

\begin{lemma} By making $C^0$ smaller, there exists a family of 
abelian surfaces
$$ a: A^0\to C^0$$
with $\rho(A^0_t)\geq 3$ such that $g:Z^0\to C^0$ is 
its Kummer construction.
\end{lemma}

\begin{proof}  It is easy to see that, for any $t\in C^0$,
$NS( Z^0_t)$ contains a sub-lattice, which is isomorphic
to 
${\Z^{19}}^{\tau}(2)\oplus\Z[N_1,...,N_8]$ 
as a trivial $\pi_1(C\setminus S)$-modules. Thus $g: Z^0\to C^0$ is a family of
Kummer surfaces with $\rho(Z_t)\geq 19$.
Let $t_0\in C^0$ with $\rho (Z^0_{t_0})=19.$ Then
$NS(Z^0_{t_0})\supset {\Z^{19}}^{\tau}\oplus\Z[N_1,...,N_8]$
and 
$$NS(Z^0_{t_0})\otimes\Q=
({\Z^{19}}^{\tau}\oplus\Z[N_1,...,N_8])\otimes\Q.$$

Let $E_1,...,E_{16}$ be the liftings of the sixteen pairwise-disjoint 
$(-2)$-curves on $Z^0_{t_0}$ to $Z^0$.
It is not difficult to see that we can choose $E_i$ ($i=1,...,16$) to
be effective divisors on $Z^0$. In fact, since
$g_*\mathcal O_{Z^0}(E_i)\not= 0$ (because $H^0(E_i|_{{Z^0}_t})\not= 0$
for any $t\in C^0$ by Riemann-Roch theorem), we have, for $m$ large 
enough and a point $p\in C^0$,
$H^0(\mathcal O_{Z^0}(E_i+m g^{-1}(p)))=
H^0(\mathcal O_{C^0}(mp)\otimes g_*\mathcal O_{Z^0}(E_i))\not=0.$
Thus there is an effective divisor $D$ on $Z^0$ such that
$D|_{Z^0_{t_0}}$ is numerical equivalent to $E_i|_{Z^0_{t_0}}$,
which implies that $D|_{Z^0_{t_0}}=E_i|_{Z^0_{t_0}}$ since
a nodal class is represented by only one effective divisor.
We can choose $E_i$ ($i=1,...,16$) to be irreducible further. In fact,
we will
show that $E_i|_{Z^0_t}$ is irreducible if $\rho(Z^0_t)=19$.
Otherwise, let $E_i|_{Z^0_t}=D_1+D_2$, where $D_1$ is irreducible
with $D_1^2=-2$ and $D_2$ is effective. Note that for any lifting of an
irreducible curve, whose restriction to any other fibre is equivalent
to an effective divisor. Thus if $\tilde D_1$ and $\tilde D_2$ are 
the liftings
of $D_1$ and $D_2$ ($\tilde D_2$ obtained by lifting the 
irreducible components
of $D_2$), we see that $\tilde D_1|_{Z^0_{t_0}}$ and
$\tilde D_2|_{Z^0_{t_0}}$ are equivalent to effective divisors.
On the other hand, $E_i|_{Z^0_{t_0}}-\tilde D_1|_{Z^0_{t_0}}$ is
numerically equivalent to $\tilde D_2|_{Z^0_{t_0}}$ since it is so on
$Z^0_t$. But this is impossible since
$E_i|_{Z^0_{t_0}}$ is a nodal class.
Let $g: Z\to C$ be a compactification of $g: Z^0\to C^0$ with
$Z$ smooth and $E_1,...,E_{16}$ be extend to $Z$. It is known
that 
$E_1|_{Z_{t_0}}+\cdots+E_{16}|_{Z_{t_0}}\equiv2\delta.$
Let $\Delta$ be a divisor on $Z$ such that
$\Delta|_{Z_{t_0}}=\delta$. Then
$E_1+\cdots +E_{16}-2\Delta$
is numerically equivalent to zero on the general fibres, thus
$$E_1+\cdots+E_{16}-2\Delta\equiv {g}^*D_a,\quad D_a\in Div(C).$$
Choose $C^0$ smaller so that $E_i|_{Z_t}$ ($i=1,...,16$)
are irreducible for $t\in C^0$ and
$$E_1+\cdots+E_{16}\equiv 2\Delta\quad \text{on $Z^0$}.$$ 
Let ${A^0}'\to Z^0$ be the double covering with branch locus
$E_1+\cdots+E_{16}$, and let $\varpi:{A^0}'\to  {A^0}$ be
the uniform blow down of the sixteen $(-1)$-curves on the fibres
${A'_0}_t$. Then $a:A^0\to C^0$ is the family of abelian surfaces 
with $\rho({A^0}_t)\geq 3$.

\end{proof}

\section{splitting on families of abelian surfaces}

Let $a: A^0\to C^0$ be the family of abelian surfaces constructed in 
Lemma 4.9. We take
a compactification $a: A\to C,$ (which may not be semi-stable). 
It is known that the sub VHS $\T_a\subset R^2a_*(\Z_{A^0})$
of the transcendental part of the weight-2
VHS attached to $a: A^0\to C^0$
is Hodge isometric to $\T_g(2),$ where
$$T_g\subset R^2g_*(\Z_{Z^0})$$
is the sub VHS of the transcendental part of the weight-2
VHS attached to $g: Z^0\to C^0.$ Furthermore, 
$T_g$ is Hodge isometric to $\T_f(2),$ where
$$\T_f=\W\subset R^2f_*(\Z_{f^{-1}(C^0)})$$
is the sub VHS of the transcendental part of the weight-2
VHS attached to $f: f^{-1}(C^0)\to C^0.$
Since $\W$ is, in fact, defined on $C\setminus S$, $\T_a$ can be extended
to $C\setminus S$ as an VHS.

\begin{lemma} The $\Q-$vector space of endomorphisms of  
$$ R^1a_*(\Z_{A^0})\otimes Q$$
has dimension 4, and is of (0,0)-type.\end{lemma}

\begin{proof}
By the construction of $a: A^0\to C^0,$ we
see $R^2a_*(\Z_{A^0})\otimes\Q $ contains a constant local system 
of dimension 3 of (1,1)-type (this corresponds to a sub-lattice of 
Picard lattice of $A^0$). Hence, it corresponds to a 3-dimensional subspace 
of $End(R^1a_*(\Z_{A^0})$ of (0,0)-type. 
Using a non-scalar endomorphism of this space, we can split 
$R^1a_*(\Z_{A^0})\otimes \C$ into the following type
$$ R^1a_*(\Z_{A^0})\otimes\C\simeq\W_1\oplus\W_2,$$ 
where both $\W_i$ are of rank-2 and irreducible over $\C.$ 
Otherwise $R^1a_*(\Z_{A^0})\otimes\C$ would contain a rank-1 
sub-local system with
zero Higgs field. This implies that the Higgs field of $(p,g)$-type on 
$\wedge^2R^1a_*(\Z_{A^0})\otimes\C$ can not be isomorphism, a contradiction.
We claim that $\W_1\simeq\W_2.$  Otherwise, the space  
$End(R^1a_*(\Z_{A^0}))\otimes \C$ has at most dimension 2, 
a contradiction. Let
$\W_1\simeq\W_2\simeq\W.$
We have then
$$\sE nd(R^1a_*(\Z_{A^0}))\otimes\C\simeq \sE nd(\W)^{\oplus 4}=
\sE nd_0(\W)^{\oplus 4}\oplus \C^4=\W'\oplus\C^4.$$
Since $ \W$ is irreducible, one shows that $\W'$ does not 
contain any constant
sub-local system and the last splitting can be defined over$\Q.$ Hence, 
$$\dim End(R^1a_*(\Z_{A^0}))\otimes\Q=4.$$ 

\end{proof}

\begin{lemma}The family $a: A^0\to C^0$ is isogenous to the square product 
of a
family of elliptic curves $e: E^0\to {C^0}.$ \end{lemma}

\begin{proof} Case 1).  Suppose that there is a subset $T\subset C^0$ of
 non-countable many points
such that $A_t$ is isogenous to $ E_t\times E_t,\, t\in T.$ Since 
there are only countable many isomorphic classes of elliptic curves 
having complex
multiplication, we find an $t_0\in T$ such that $End(E_{t_0})\otimes Q=\Q.$ 
Hence, the endomorphism algebra
$$ End(A_{t_0})\otimes Q\simeq M_2(\Q).$$
In the other words, we have
$ End(R^1a_*(\Z_{A^0})\otimes\Q)|_{t_0}\simeq M_2(\Q).$
Since 
$$End(R^1a_*(\Z_{A^0})\otimes\Q)$$ is constant local system, we have
$ End(R^1a_*(\Z_{A^0})\otimes\Q)\simeq M_2(\Q).$
The element 
$$ \left[\begin{array}{cccccccccc}
1& 0\\
0&0\\
\end{array}\right]\in End(R^1a_*(\Z_{A^0})\otimes\Q) $$ 
gives a $\Q-$splitting  
$ R^1a_*(\Z_{A^0})\otimes\Q=\W_{\Q}\oplus\W_{\Q},$
thus isogeny splitting of $f: A^0\to C^0$
into the square product of a family of elliptic curves
$ e: E^0\to {C^0}.$\\

Case 2). Suppose that there are non-countable many points
 $\{t\}\subset C^0$ such that $A_t$ is simple. Since
the Picard number $\rho(A_t)\geq 3,$ one checks easily that
$\rho(A_t)=3$ and $End(A_t)\otimes\Q$ is the totally indefinite 
quaternion algebra
over $\Q.$ An abelian surface with this type endomorphism algebra is called
a false elliptic curve. There are countable many projective curves 
$\{C_i\}_{i\in\N}$ in the moduli space
of polarized abelian surfaces, which are Shimura curves of certain type 
and parametrice all false elliptic curves. So, the family
$a: A^0\to C^0$ induces a morphism $\phi: C^0\to C_i$ for some $i\in \N,$
which extends to a morphism $\phi: C\to C_i.$ This implies that 
the local monodromies of $R^2a_*(\Z_{A^0})$ around the singularity 
has finite order.
It contradicts to $S\not=\emptyset.$

\end{proof}

\section{Proof of theorems and corollary}

{\bf Proof of Theorem 0.1}:$\quad$ Only the modularity of 
$C'\setminus \sigma^{-1}S$ needs to be checked. 
The isogeny 
$a: A^0\to C^0\sim e^2: E^0\times_{C^0} E^0\to C^0$
induces an isomorphism
$ S^2(R^1e_*(\Z_{E^0}))\simeq\W|_{C^0}$.
There are natural group homomorphisms
$$ 1\to\{\pm 1\}\to SL_2(\R)\to SO(1,2),$$
which induce an isomorphism between $\sH$ and a connected
component of the symmetric space $SO(1,2)/SO(2)\times O(1),$ say
$$ i:\sH\simeq SO^+(1,2)/SO(2)\times O(1).$$

Since $\W|_{C^0}$ is the restriction of $\W$ on $C\setminus S$ to
$C^0,$ the local monodromies of $R^1e_*(\Z_{E^0})$ around 
$(C\setminus S)\setminus C^0$ are either $+1,$ or $-1$. Thus the 
projective 
monodromy representation of $R^1e_*(\Z_{E^0})$ is actually
defined on $C\setminus S,$ say
$$\rho_{R^1e_*\Z_{E^0}}:\pi_1(C\setminus S,*)\to PSL_2(\Z).$$ 
Let 
$\widetilde \phi_{R^1e_*(\Z_{E^0})}:\widetilde{C\setminus S}\to \sH$
be the period map corresponding to $R^1e_*\Z_{E^0}$ and 
$$\widetilde \phi_{\W}:\widetilde {C\setminus S} 
\to SO^+(1,2)/SO(2)\times O(1))$$  
denote the period map corresponding to $\W$. Then 
$\widetilde\phi_{\W}=i\cdot\widetilde\phi_{R^1e_*(\Z_{E^0})}$ is
an isomorphism. In fact, the tangent map of $\widetilde \phi_{\W}$
is precisely the Kodaira-Spencer map of $\W$:
$\theta^{2,0}: E^{2,0}\to E^{1,1}_1\otimes \Omega^1_C(\log S),$
which is isomorphic at each point by Lemma 1.1. Thus $\widetilde \phi_{\W}$
is a local diffeomorphism. Since
the Hodge metric on the Higgs bundle corresponding to $\W$ has 
logarithmic growth at $S$ and bounded curvature by Schmid \cite{Sch}, 
together with the remarks after Proposition 9.1 and Proposition 9.8
in \cite{Si3}, $\widetilde \phi_{\W}$ is a covering map, 
hence an isomorphism.
This implies that  $\widetilde \phi_{R^1e_*(\Z_{E^0})}$ is an isomorphism. 
Thus
$$ \phi_{R^1e_*(\Z_{E^0})} : 
C\setminus S\simeq \sH/ \rho_{R^1e_*(\Z_{E^0})}$$
is an isomorphism.\\

In order to prove Theorem 0.2, we need the following lemma.

\begin{lemma} Let $f:X\to C$ be a family of semi-stable K3 surfaces, which has
zero iterated Kodaira-Spencer map and reaches
the Arakelov bound (II)
$$\deg f_*\omega_{X/C}=\frac{1}{2}\deg\Omega^1_C(\log S).$$

Then, after passing through a finite {\'e}tale covering $C'\to C$, the VHS
$\W$ is non-rigid.
\end{lemma}

\begin{proof} One needs to show that,
after passing through a finite {\'e}tale covering of $C$, the local
system $R^2f_*(\Z_{X^0})\otimes\C$ admits a non-zero endomorphism of 
type $(-1, 1).$
By Lemma 3.2, one has splitting
$$R^2f_*(\Z_{X^0})\supset\W_{\Z}\oplus\Z^{18},
\quad R^2f_*(\Z_{X^0})\otimes\Q
=(\W_{\Z}\oplus\Z^{18})\otimes\Q.$$

By Lemma 1.2, the Higgs bundle corresponds to $\W$ has the form
$$ (E^{2,0}\oplus E^{1,1}_1)\oplus({E^{1,1}_1}^*\oplus E^{0,2})$$
such that the Higgs fields
$$ \tau: E^{2,0}\to E^{1,1}_1\otimes\Omega^1_C(\log S),
\quad \tau^*: {E^{1,1}}^*\to E^{0,2}\otimes\Omega^1_C(\log S)$$
are isomorphisms. These two Higgs subbundles correspond to
two sub-local systems $\W_1$ and $\bar W_1$. 
We claim that, after passing through a finite {\'e}tale covering 
of $C$, one has $\W_1\simeq \bar\W_1$. To prove the claim, consider the 
sub-local system $$ \W_1\to \W.$$
If $\W_1$ is not rigid, then there is a small deformation 
$W_{1,t}\subset\W\otimes\C$ 
such that both projections
$ W_{1,t}\subset\W\otimes\C\to \W_1$ and
$W_{1,t}\subset\W\otimes\C\to \bar\W_1 $
are non-zero. Since $\W_1$ is irreducible, one obtains 
$$ \W_1\simeq\W_{1,t}\simeq\bar\W_1.$$

If $\W_1$ is rigid, then by Lemma 2.1 $\W_1$ is defined over 
a number field 
$K$. Let $\mathcal O_K$ denote the ring of algebraic integers in $K$, and
let 
$$\W_{1\mathcal O_K}=\W\otimes_{\Z}\mathcal O_K\cap\W_1.$$

Then $\W_{1\mathcal O_K}\otimes K=\W_1$,  
which means that the corresponding monodromy representation of $\W_1$ 
can be defined
over $\mathcal O_K$. The determinant $\det\W_1=E^{2,0}\otimes E^{1,1}_1$ 
is a rank-1 unitary local system $\eta\in Pic^0(C)$ and takes values 
in $\mathcal O_K$. By a theorem of Kronecker, $\eta$ is a torsion. 
So, after passing through the finite {\'e}tale covering corresponding to
$\eta$, one obtains $E^{2,0}\simeq {E^{1,1}_1}^*$ and
$$(E^{2,0}\oplus E^{1,1}_1,\tau)\simeq (E^{1,1*}_1\oplus E^{0,2},\tau^*).$$

Thus, in any case, we obtain a non-zero endomorphism
$$(E^{2,0}\oplus E^{1,1}_1)\oplus({E^{1,1}_1}^*\oplus E^{0,2})
\to (E^{2,0}\oplus E^{1,1}_1)\oplus({E^{1,1}_1}^*\oplus E^{0,2})$$

of type $(-1,1)$, which corresponds to an endomorphism of 
$R^2f_*(\Z_{X^0})\otimes\C$ of
type $(-1,1)$. 

\end{proof}

{\bf Proof of Theorem 0.2}:$\quad$ By Lemma 6.1, after passing through a 
finite
{\'e}tale covering $C'\to C$, the VHS $\W$ is non-rigid. 
By Corollary 5.6.3 of \cite{S-Z}, one has
$$ End(\W)\otimes\Q\simeq M_2(\Q).$$

Taking an element in $M_2(\Q)$ with two distinct rational eigenvalues, 
we get a $\Q-$splitting
$\W\otimes\Q=\W_1\oplus\W_2$
such that $\W_1$ is isomorphic to $\W_2$ and the Higgs bundle 
corresponding to $\W_1$ has the form
$$ (L\oplus L^{-1},\theta),\quad \theta: L\simeq L^{-1}\otimes
\Omega^1_C(\log S).$$ 

$\W_1$ has an $\Z-$structure defined by
$ \W_{1\Z}=\W_{\Z}\cap \W_1$.
Again by Proposition 9.1 of \cite{Si3}, the Higgs bundle 
$\theta: L\simeq L^{-1}\otimes\Omega^1_C(\log S) $ 
gives rise to the uniformization
$$ C\setminus S\simeq \sH/\rho_{\W_1}\pi_1(C\setminus S,*),$$

where $\rho_{\W_1}\pi_1(C\setminus S,*)\subset SL_2(\Z)$ of finite index.\\

{\bf Proof of Corollary 0.4} i) 
By Theorem 0.1 there exists a family
of elliptic curves 
$g: E^0\to {\P^1}^0\subset \P^1\setminus S$
such that the projective representation 
$$ p\rho_{R^1g_*{\Z_{E^0}}}:\pi_1(\P^1\setminus S,*)\to 
\Gamma'\subset PSL_2(\Z)$$
extends to $\P^1\setminus S$
and $ \P^1\setminus S\simeq \sH/\Gamma'$. 
By \cite{B2}, $\Gamma'\subset PSL_2(\Z)$ is of index 12 and conjugates 
to one of the following 6 subgroups of $PSL_2(\Z),$ which are
images of $\Gamma(3)$, $\Gamma^0_0(4)\cap\Gamma(2)$, $\Gamma^0_0(5)$,
$\Gamma^0_0(6)$, $\Gamma_0(8)\cap\Gamma_0^0(4)$ and 
$\Gamma_0(9)\cap\Gamma^0_0(3)$ in $SL_2(\Z)$
of index 24, where
$$ \Gamma(n)= \left\{  \left[ \begin{array}{ccccccccc}
a & b\\
c & d\\
\end{array}\right]\in SL_2(\Z)|b\equiv c\equiv 0,\, a\equiv 1 (\text {mod.}n)\right\}, $$

$$\Gamma^0_0(n)=\left\{  \left[ \begin{array}{ccccccccc}
a & b\\
c & d\\
\end{array}\right]\in SL_2(\Z)| c\equiv 0,\, a\equiv 1 (\text {mod.}n)\right\}, $$

$$ \Gamma_0(n)=\left\{  \left[ \begin{array}{ccccccccc}
a & b\\
c & d\\
\end{array}\right]\in SL_2(\Z)| c\equiv 0 (\text {mod.}n)\right\}. $$  

In the proof of Theorem 0.1, we have seen already that the  monodromy of 
$R^1g_*\Z_{E^0}$ of
a short loop around a point of $(\P^1\setminus S)\setminus {\P^1}^0$ is 
either
$+1$, or $-1$. If all of them equal to $+1$, then the representation
$\rho_{R^1g_*\Z_{E^0}}$ extends to $\P^1\setminus S$, and the image of 
$\pi_1(\P^1\setminus S,*)$
under this representation conjugates to one of the above 6 subgroups. 
Hence $g: E^0\to {\P^1}^0$ extends to 
a modular family of elliptic curves $g: E\to \P^1\setminus S$ from one 
of 6 examples in \cite{B2}.
Suppose that the  monodromies  of $R^1g_*\Z_{E^0}$ of
 short loops around some points of 
$(\P^1\setminus S)\setminus {P^1}^0$ equal to
$-1$. Then the image of $\pi_1(\P^1\setminus S,*)$ conjugates 
to the preimage
$p^{-1}p\Gamma$, where $\Gamma$ is one of 
 $\Gamma(3)$, $\Gamma^0_0(4)\cap\Gamma(2)$, $\Gamma^0_0(5)$,
$\Gamma^0_0(6)$, $\Gamma_0(8)\cap\Gamma_0^0(4)$ and 
$\Gamma_0(9)\cap\Gamma^0_0(3)$.
The inclusion $\Gamma\subset p^{-1}p\Gamma$ of index 2 defines an 
{\'e}tale covering
${E^0}'\to E^0,$ which is {\'e}tale along the fibres and the family 
$g': {E^0}'\to
{\P^1}^0$ extends to the modular family of elliptic curves 
$g': E'\to \P^1\setminus S$
corresponding to $\Gamma$.\\
ii) is straightforward.

\bibliographystyle{plain}

\begin{thebibliography}{XXX} % XXX=breitestes Label
% \bibitem{} gibt [1] \bibitem[XX]{} gibt [XX]
\bibitem{B1} Beauville, A.:
Le nombre minimum de fibres singulier\`es d'une courbe stable sur
$P\sp 1.$ (French) Asterisque {\bf 86} (1981) 97-108
\bibitem{B2} Beauville, A.:
Les familles stables de courbes elliptiques sur $P\sp{1}$
admettant quatre fibres singuli\`eres. C. R. Acad. Sci. Paris Sér.
I Math. {\bf 294} (1982) 657--660.
\bibitem{BL}  Birkenhage, Ch., Lange, H.: Complex Abelian Varieties.
Grundlehren d. math. Wiss. {\bf 302}. Springer-Verlag, Berlin-Heidelberg
1992
\bibitem{Del} Deligne, P.: Th\'eorie de Hodge II. I.H.\'E.S. Publ. Math.
{\bf 40} (1971) 5--57.
\bibitem{Do} Dolgachev, I.: Mirror symmetry for lattice polarized $K3$ surfaces. Preprint,
AG/9502005.
\bibitem{F} Fujita, T. : On K\"ahler fibre spaces over curves.
 J. Math. Soc. Japan. {\bf 30} (1978) 779--794.
\bibitem{G}Griffiths, P.: Topics in transcendental algebraic geometry.
Ann of Math. Stud. 106 Princeton Univ. Press. Princeton, N.J. 1984.
\bibitem{J-Z} Jost, J, Zuo, K.: Arakelov type inequalities for
Hodge bundles over algebraic varieties, Part 1: Hodge bundles
over algebraic curves. preprint (1999), to appear in J. Alg. Geometry
\bibitem{K}Kawamata, Y.: Characterization of abelian varieties. Comp. Math. {\bf 43} (1981) 253--276.
\bibitem{Li-Y1} Lian, B, Yau, S-T.:  Arithmetic properties of mirror map and quantum coupling. Preprint 1994, hep-th.
\bibitem{Li-Y2} Lian, B, Yau, S-T.:  Mirror maps, modular relations and hypergeometric
series II, Preprint 1994, hep-th.
\bibitem{L-P} Looijenga, E, Peters, C.: Torelli theorems
for K{\"a}hler K3 surfaces. Comp. Math.{\bf 42}1981) 145--186
\bibitem{Mo1} Mok, N.: Metric rigidity theorems on Hermitian locally symmetric manifolds.
Series in Pure Mathematics Vol 6, World Scientific, Singapore-New Jersey-London-Hong Kong, 1989.
\bibitem{Mo2} Mok, N.: Aspects of K\:ahler geometry on arithmetic varieties. A.M.S. Proc.
of Symp. in Pure Math. Vol 52. (1991), Part 2, 335-396.
\bibitem{Mo} Morrison, D:
On K3 surfaces with large Picard number. Invent. Math.{\bf 75} (1984) 105--121.
\bibitem{Mum} Mumford, D.: Abelian Varieties. (1970) Oxford Univ. Press,
Oxford
\bibitem{N} Nikulin, V.: Finite groups of automorphisms of K\"ahler surfaces of type $K3$.
Trudy Mosk. Math. Ob. {\bf 38} (1979) 75-137. Trans. Moscow Math. Soc. {\bf 38} (1980) 71-135.
\bibitem{Pet} Peters, C.: Arakelov-type inequalities for Hodge
bundles. preprint (1999)
\bibitem{Sai} Saito, M.-H.: Classification of nonrigid
families of abelian varieties.  Tohoku Math. J. {\bf 45}  (1993)
159--189.
\bibitem{S-Y} Saito, M-H., Yui, N.: The modularity conjecture for rigid Calabi-Yau
three-fold over $\Q.$ Preprint, AG/0009041.
\bibitem{S-Z} Saito, M.-H., Zucker, S.:
Classification of nonrigid families of $K3$ surfaces and a
finiteness theorem of Arakelov type.  Math. Ann.  {\bf 289} (1991)
1--31
\bibitem{Sch} Schmid, W.: Variation of Hodge structure: The singularities
of the period mapping. Invent. math. {\bf 22} (1973) 211--319
\bibitem{Sim} Simpson, C.: Higgs bundles and local systems. Publ. Math.
I.H.E.S {\bf 75} (1992) 5--95
\bibitem{Si2} Simpson, C.: Harmonic bundles on noncompact curves.
Journal of the AMS {\bf 3} (1990) 713--770
\bibitem{Si3} Simpson, C.: Constructing variations of Hodge structure using Yang-Mills theory
and applications to uniformization. Journal of the AMS {\bf 1} (1988) 867-918.
\bibitem{Sha} Shafarevich, I. R.: On some families of abelian surfaces.
Izv. Math. {\bf 60} (1996) 1083--1093
\bibitem{Shi} Shimura, G.: Introduction to the Arithmetic Theory
of Automorphic Functions. Publ. Math. Soc. of Japan 11 (1971),
Iwanami Shoten and Princeton University Press
\bibitem{Tan} Tan, S-L.: The minimal number of singular fibers of a
semistable curve over $P\sp 1$.
J. Alg. Geom. {\bf 4} (1995) 591--596
 \bibitem{V} Viehweg, E: Weak positivity and additivity of Kodaira
dimension for certain algebraic fibre spaces. Adv. Studies in Pure Math. 1, Algebraic
Varieties and Analytic Varieties (1983) 329--353
\bibitem{V-Z} Viehweg, E., Zuo K.:  Families of abelian varieties over curves with maximal Higgs field. Preprint math.AG/0204261.
\bibitem{Y}Yau, S-T.: A general Schwarz lemma for Kahler manifolds Amer. J. Math.
{\bf 100} (1978) 197-203.
\end{thebibliography}

\enddocument